\documentclass[noamsfonts]{amsproc}
\usepackage{aspmproc}
\title[Logarithmic forms and anti-invariant forms of
reflection groups]{Logarithmic forms and anti-invariant forms of
reflection groups}
\author[A. Shepler and H. Terao]{Anne Shepler
and
Hiroaki Terao$^{1}$\footnote{$^{1}$
Partially supported by
NSF-DMS 9504457.}
\\
{\footnotesize 
Dedicated to Peter Orlik on his sixtieth birthday
}}

\address{{\rm Anne Shepler}\\ 
Mathematics Department\\
University of California at San Diego\\
La Jolla, CA 92093 USA\\
{\rm Hiroaki Terao}\\
Mathematics Department\\
Tokyo Metropolitan University\\
Hachioji, Tokyo 192-0372 Japan\\
}

\email{
ashepler@euclid.ucsd.edu hterao@comp.matro-u.ac.jp
}
 
\theoremstyle{plain} 
\newtheorem{theorem}{Theorem}

\newtheorem{proposition}[theorem]{Proposition}
\newtheorem{corollary}[theorem]{Corollary}

\theoremstyle{definition}

\newtheorem{example}[theorem]{Example}

\theoremstyle{remark}

\newcommand{\A}{{\mathcal A}}

\begin{document}

\begin{abstract}
Let $W$ be a finite group generated by unitary
reflections and $\A$ be the set
of reflecting hyperplanes.  We will give a characterization of the
logarithmic differential forms with poles along
$\A$ in terms of anti-invariant differential forms. 
If $W$  is a Coxeter group defined over ${\mathbf R}$, then
the characterization provides a new method to
find a basis for the module of logarithmic 
differential forms out of basic invariants. 
\end{abstract}
\maketitle

{\bf Basic definitions.}
Let $V$ be an $\ell$-dimensional unitary space.
Let $W \subset {\bf GL}(V)$ be a finite group generated by unitary
reflections and $\A$ be the set
of reflecting hyperplanes.  We say that $W$ is  a {\em unitary reflection 
group}
 and $\A$
is the corresponding {\em unitary reflection arrangement}.
Let $S$ be the
algebra of polynomial functions on $V.$ The algebra $S$ is naturally
graded by $S = \bigoplus_{q \geq 0}~ S_q$ where $S_q$ is the space of
homogeneous polynomials of degree $q.$ Thus $S_1 = V^*$ is the dual
space of $V.$ Let ${\rm Der}_S$ be the $S$-module of ${\bf
C}$-derivations of $S.$ We say that $\theta \in {\rm Der}_S$ is
homogeneous of degree $q$ if $\theta (S_1) \subseteq S_q.$ Choose for each
hyperplane $H \in {\mathcal A}$ a linear form $\alpha_H \in V^*$ such that
$H = {\rm ker}(\alpha_H).$  Define $Q \in S$ by
\[
       Q = \prod_{H \in {\mathcal A}}~ \alpha_H \, .
\]
The polynomial $Q$ is uniquely determined, up to a constant
multiple, by the group $W.$  
When convenient we choose a basis $e_1,\ldots,e_l$ for
$V$ and let $x_1,\ldots,x_l$ denote the dual basis for $V^*.$  Let
$\langle ~~, ~~\rangle : V^* \times V \rightarrow {\bf C}$ denote the natural
pairing. Thus $\langle x_i, e_j \rangle = \delta_{ij}.$ 
For each $v\in V$ let $\partial_v \in {\rm Der}_S$ 
be the unique derivation such
that $\partial_v x = \langle x, v \rangle$ for $x \in V^*.$ Define
$\partial_i \in {\rm Der}_S$ by $\partial_i = \partial_{e_i}.$ Then
$\partial_i x_j = \delta_{ij}$ and ${\rm Der}_S$ is a free $S$-module
with basis $\partial_1,\ldots,\partial_l.$   There is a natural
isomorphism $S \otimes V \rightarrow {\rm Der}_S$ of $S$-modules given
by
\[
          f \otimes v \mapsto f \partial_v
\]
for $f \in S$ and $v \in V.$  Let $\Omega^1 = {\rm
Hom}_S({\rm Der}_S,S)$ be the $S$-module dual to ${\rm Der}_S.$ Define
$d: S \rightarrow \Omega^1$ by $df(\theta) = \theta(f)$ for $f \in S$
and $\theta \in {\rm Der}_S.$ Then $d(ff') = (df)f' + f(df')$ for $f,f'
\in S.$ Furthermore, $\Omega^1$ is a free $S$-module with basis
$dx_1,\ldots,dx_l$ and $df = \sum_{i=1}^l~ (\partial_if) dx_i.$ There is
a natural isomorphism $S \otimes V^* \rightarrow \Omega^1$ of
$S$-modules given by
\[
          f \otimes x \mapsto f dx
\]
for $f \in S$ and $x \in V^*.$  The modules ${\rm
Der}_S$ and $\Omega^1$ inherit gradings from $S$ which are defined by
${\rm deg}(f \partial_v ) = {\rm deg}(f)$ and ${\rm deg}(f dx) = {\rm
deg}(f)$ if $f \in S$ is homogeneous.
Define $\Omega^{p} = \bigwedge^{p}_{S} \Omega^{1}  ~~
(p = 1, \ldots ,\ell)  $. 
Let $\Omega^{0} = S$. 
The $S$-module $\Omega^{p} $ is free with a basis 
$\{dx_{i_{1} }\wedge\cdots\wedge dx_{i_{p} }  \mid 1\leq i_{1} <\cdots < i_{p}\leq \ell  \}$.
It is naturally isomorphic to 
$S\otimes_{{\bf C}} \bigwedge^{p} V^{*}    $.  
   Let $\Omega^p(\A)$ be the $S$-module
of {\em logarithmic $p$-forms} with poles along $\A$ \cite{sai3}\cite{ort}:
\[
\Omega^{p}(\A) = \{ \frac{\eta}{Q} \mid \eta\in \Omega^{p}, 
d(\frac{\eta}{Q}) \in \frac{1}{Q}  \Omega^{p+1}\} 
\]
where $d$ is the exterior differentiation.

The unitary reflection group $W$ acts contragradiently on $V^{*} $
and thus on $S$.  The modules 
${\rm Der}_{S}$ and $\Omega^{p}~~(p = 0, \ldots ,\ell)$
also have $W$-module structures so that the above isomorphisms
are $W$-module isomorphisms.  
If $M$ is an ${\bf C}[W]$-module let
$M^W = \{ x \in M ~|~ wx = x ~{\rm for~all}~ w \in W \}$
denote the space of invariant elements in $M.$  Let
$M^{\det^{-1} } = \{ x \in M ~|~ wx = \det(w)^{-1} x ~{\rm for~all}~ w \in W \}$
denote the space of anti-invariant elements in $M.$  
Let $R = S^{W} $ be the invariant subring of $S$ under $W$.
By a theorem of Shephard,
Todd, and Chevalley \cite[V.5.3, Theorem 3]{bou} there exist algebraically
independent homogeneous polynomials $f_1,\ldots,f_l \in R$ such that $R
= {\bf C}[f_1,\ldots,f_l].$  They are called {\em basic invariants}.
Elements of $S^{\det^{-1} } $ and $(\Omega^{p})^{\det^{-1} }  $ are called 
{\em anti-invariants} and {\em anti-invariant
$p$-forms} respectively.   
It is well-known that $S^{\det^{-1} } = RQ $.

\bigskip

{\bf The main theorem.}
The following theorem gives the relationship between
logarithmic forms and anti-invariant forms.   

\medskip
\begin{theorem}
\label{thm1}
For $0\leq p\leq\ell$, 
$$\Omega^p(\A) = \frac{1}{Q}(\Omega^p)^{\det^{-1} } \otimes_{R} S.$$
\end{theorem}
 
\begin{proof} 
When $p=0$, the result follows from the formula
$S^{\det^{-1} } = RQ $.
 Let $p > 0$. 
Let $x_1, \ldots , x_\ell$ be an orthonormal basis
for $V^*$.  Let $\theta_1, \ldots, \theta_{\ell}$ be an $R$-basis for
${\rm Der}_{S}^{W} $.
Then, by \cite[Theorem 6.59]{ort},
 $\theta_1, \ldots, \theta_{\ell}$ is known to be an $S$-basis
 for the module $D(\A)$ of $\A$-derivations,
 where
 \[
D(\A) = \{\theta\in {\rm Der_{S}} \mid \theta(Q) \in QS\}.
 \]
 By the contraction  
 $\left< ~, ~\right>$ of a $1$-form and a derivation,
 the $S$-modules
 $D(\A)$ and $\Omega^1 (\A)$ are $S$-dual
each other \cite[p.268]{sai3} \cite[Theorem 4.75]{ort} .   Let
$\{\omega_{1}, \ldots , \omega_{\ell}  \}
\subset \Omega^{1}(\A) $
be dual to  
$\{\theta_{1}, \ldots , \theta_{\ell}  \}$.
In other words, 
$\left<\theta_{i}, \omega_{j}  \right> =
\delta_{ij} $ (Kronecker's delta).
Then  
$\{\omega_{1}, \ldots , \omega_{\ell}  \}
$
is an $S$-basis for $\Omega^{1}(\A) $.
Then each $\omega_{i} $ is obviously $W$-invariant
and
\[
\omega_{i} \in (\frac{1}{Q} \Omega^{1}   )^{W}
=
 \frac{1}{Q} (\Omega^{1} )^{\det^{-1} }. 
\]
   Therefore we have
   \[
   \Omega^{1}(\A) \subseteq \frac{1}{Q} (\Omega^{1} )^{\det^{-1} } \otimes_{R} S.    
   \]
 By \cite[Proposition 4.81]{ort}, the set
 $\{\omega_{i_{1} }\wedge\cdots\wedge\omega_{i_{p} } \mid
 1 \leq i_{1} < \cdots < i_{p} \leq \ell \}$
 is a basis for $\Omega^{p}(\A) $.
 In particular,    
$   \omega_{i_{1} }\wedge\cdots\wedge\omega_{i_{p} }
\in \frac{1}{Q} \Omega^{p}   $.
Since $\omega_{i_{1} }\wedge\cdots\wedge\omega_{i_{p} }$
is $W$-invariant,
$Q(\omega_{i_{1} }\wedge\cdots\wedge\omega_{i_{p} })
\in
(\Omega^{p} )^{\det^{-1} } $.
This shows that
\[
\Omega^{p}(\A) \subseteq \frac{1}{Q} (\Omega^{p} )^{\det^{-1} }\otimes_{R} S.      
\]
   
   Conversely let $\omega\in \frac{1}{Q} (\Omega^{p} )^{\det^{-1} }  $.
   Then $Q\omega \in \Omega^{p} \subseteq \Omega^{p}(\A) $.
   Thus $Q\omega$ can be uniquely expressed as
   \[
   Q\omega = 
   \sum_{i_{1} < \cdots < i_{p} } f_{i_{1}\cdots i_{p}  }
    \omega_{i_{1} }\wedge\cdots\wedge\omega_{i_{p} }
    ~~~(f_{i_{1}\cdots i_{p}  } \in S).
   \]
   Act $w\in W$ on both sides, and we get
   \[
   \det(w)^{-1} Q\omega = w(Q) \omega = 
\sum_{i_{1} < \cdots < i_{p} } w(f_{i_{1}\cdots i_{p}  })
    \omega_{i_{1} }\wedge\cdots\wedge\omega_{i_{p} }.
        \]
        Therefore, by the uniqueness of the expression, 
        we have 
        \[
        \det(w)^{-1}  f_{i_{1}\cdots i_{p}  }
        =
        w(f_{i_{1}\cdots i_{p}  })~~~~~(w\in W)
        \]
        and $f_{i_{1}\cdots i_{p}  } \in S^{\det^{-1} } = R Q$.
        This implies that each $f_{i_{1}\cdots i_{p}  }/Q$ 
        lies in $S$ and that 
        \[
        \omega = 
   \sum_{i_{1} < \cdots < i_{p} } 
   \left(\frac{f_{i_{1}\cdots i_{p}  }}{Q}\right) 
    \omega_{i_{1} }\wedge\cdots\wedge\omega_{i_{p} } \in
    \Omega^{p}(\A).
        \]
        Thus we have shown the inclusion
         \[
         \frac{1}{Q}(\Omega^{p} )^{\det^{-1} } \otimes_{R} S
         \subseteq 
         \Omega^{p}(\A).
         \]
\end{proof}

Taking the $W$-invariant parts of the both sides in Theorem~\ref{thm1},
we have 

\begin{corollary}
\label{cor2} 
For $0\leq p\leq\ell$, 
\[
(\Omega^{p}(\A))^{W}  
=
\frac{1}{Q}  (\Omega^{p} )^{\det^{-1}  } .
\]
\end{corollary}

The following theorem is a special case of a theorem obtained
by Shepler \cite{she1}.

\begin{theorem}[Shepler]
\label{thm3} 
  For $0\leq p\leq\ell$, 
\[
(\Omega^{p} )^{\det^{-1} } = Q^{1-p} \bigwedge_{R}^{p} 
(\Omega^{1} )^{\det^{-1} } .    
\]
\end{theorem}
 
\begin{proof} Let $p=0$.
We naturally interpret the ``empty exterior product''
to be equal to the coefficient ring.
Thus the result follows from the formula
$S^{\det^{-1} } = RQ $.
Let $p > 0$.
   In the proof of Theorem~\ref{thm1}, we have already shown that
the both sides have the same $R$-basis
\[
\{Q(\omega_{i_{1} }\wedge\cdots\wedge\omega_{i_{p} }) \mid
 1 \leq i_{1} < \cdots < i_{p} \leq \ell \}.   
\]
\end{proof} 

\medskip

{\bf The Coxeter case.} From now on we assume that
$W$ is a {\em Coxeter group}.
In other words, for an $\ell$-dimensional
Euclidean space $V$, $W\subset {\bf GL}(V)$
is a finite group generated by orthogonal reflections
and  $W$ acts irreducibly on $V$.
The objects like $S$, $R$, and $\Omega^{p} $ are defined
over ${\bf R} $.
Note that $\det(w)$ is either $+1$ or $-1$ for any $w\in W$
and thus $\det = \det^{-1} $.
  
Recall the definition of the ${\bf R}$-linear map
$
\hat{d} : S \longrightarrow \Omega^{1} 
$ in \cite[Proposition 6.1]{sot}:
\[
\hat{d}f = \sum_{i=1}^{\ell} (\partial_{i} f ) d(Q(Dx_{i} )).  
\]
Here $D$ is a Saito derivation introduced in
\cite{sai2}\cite{sys}.    
The following proposition is Proposition 6.1 in \cite{sot}:

\begin{proposition}[Solomon-Terao]
\label{prop4} 
Let $f_{1} , \ldots , f_{\ell} $ be basic invariants.
Then
\[
(\Omega^{1})^{\det} = R \hat{d}f_{1} \oplus \cdots \oplus R \hat{d}f_{\ell}.    
\]
\end{proposition}
 
 From Theorem~\ref{thm3}  and Proposition~\ref{prop4}  we get

\begin{corollary}
\label{cor5} 
For $0\leq p\leq\ell$, 
\[
(\Omega^{p})^{\det} = \bigoplus_{1\leq i_{1} < \cdots < i_{p} \leq \ell} 
R Q^{1-p}  (\hat{d}f_{i_{1} }\wedge\cdots\wedge \hat{d} f_{i_{p} }   ).
\]
\end{corollary}

Using Theorem~\ref{thm1}, we have

\begin{corollary}
\label{cor6} 
For $0\leq p\leq\ell$, 
\[
\Omega^{p}(\A) = \bigoplus_{1\leq i_{1} < \cdots < i_{p} \leq \ell} 
SQ^{-p}  (\hat{d}f_{i_{1} }\wedge\cdots\wedge \hat{d} f_{i_{p} }   ).
\]
\end{corollary}

This corollary gives a new method using the new differential
operator $\hat{d} $ to calculate a basis for 
the module of logarithmic forms.

Taking the $W$-invariant parts of the both sides in Corollary~\ref{cor6},
we also have 

\begin{corollary}
For $0\leq p\leq\ell$, 
\[
(\Omega^{p}(\A))^{W}  = \bigoplus_{1\leq i_{1} < \cdots < i_{p} \leq \ell} 
R Q^{-p}  (\hat{d}f_{i_{1} }\wedge\cdots\wedge \hat{d} f_{i_{p} }   )
.
\]
\end{corollary}

\begin{example}[$B_{2}$]
When $W$ is the Coxeter group of type $B_{2} $, we can choose
\[
f_{1} = \frac{1}{2} (x_{1}^{2} + x_{2}^{2}),~~~~
f_{2} = \frac{1}{4} (x_{1}^{4} + x_{2}^{4}).  
\]
  Then, as seen in \cite[\S 5.2]{sot}, the operator $\hat{d}$
in Proposition~\ref{prop4}  satisfies
\[
\hat{d}x_{1} = -dx_{2},~~~~
\hat{d}x_{2} =  dx_{1}.  
\]
Thus
\[
\hat{d}f_{1} = -x_{1} dx_{2} + x_{2} dx_{1},~~~~
\hat{d}f_{2} = -x_{1}^{3}  dx_{2} + x_{2}^{3}  dx_{1}.  
\]
Then $\hat{d}f_{1}  $ and $\hat{d}f_{2}   $   form an $R$-basis for
$(\Omega^{1} )^{\det} $ 
and
  $\hat{d}f_{1}/Q  $ and $\hat{d}f_{2}/Q   $   form an $S$-basis for
$\Omega^{1}(\A) $  
as Corollaries~\ref{cor5}  and \ref{cor6}  assert.
\end{example}

\end{document}